\documentstyle[12pt]{article}
\begin{document}

\newcommand{\mysection }[1]{\section{#1}\setcounter{equation}{0}}

\title{\bf An Extension of Mok's Theorem on the Generalized Frankel Conjecture} \vspace{15mm}
\author{\bf}
\date{}
\maketitle \centerline{\large \bf Hui-Ling Gu and Zhu-Hong Zhang }
\vspace{8mm} \centerline{{\large Department of Mathematics}}
\vspace{3mm} \centerline{{\large Sun Yat-Sen University }}
\vspace{3mm} \centerline{{\large Guangzhou, P.R.China}}
\vspace{3mm}
%\centerline{{revised \today}}
\vspace{15mm}

\noindent {\bf Abstract } In this paper, we will give an extension
of Mok's theorem on the generalized Frankel conjecture under the
condition of the orthogonal holomorphic bisectional curvature.

\begin{center}

{\bf \large \bf{1. Introduction}}
\end{center}

Let $M^n$ be a complex $n$-dimensional compact K$\ddot{a}$hler
manifold. One of the interesting problems is to give the
classification of the manifolds under certain curvature conditions.
Corresponding to the sectional curvature condition in Riemannian
geometry, one usually considers the holomorphic bisectional
curvature in complex differential geometry. In 1979 Mori \cite{Mori}
and in 1980 Siu-Yau \cite{Siu-Yau} independently proved the famous
Frankel conjecture by using different methods. They proved that:
\emph{any compact K$\ddot{a}$hler manifold with positive holomorphic
bisectional curvature must be biholomorphic to the complex
projective space.} After the work of Mori and Siu-Yau, in 1988, Mok
\cite{Mok} generalized the Frankel conjecture to the nonnegative
case, usually we call it the generalized Frankel conjecture which
states that: \emph{any compact irreducible K$\ddot{a}$hler manifold
with nonnegative bisectional curvature must be either a Hermitian
symmetric manifold or biholomorphic to the complex projective
space.} Recently, based on the work of Brendle-Schoen \cite{BS2},
the first author \cite{Gu} gave a simple and completely
transcendental proof to Mok's theorem on the generalized Frankel
conjecture. In the late 80's, Cao and Hamilton \cite{CH} introduced
the concept of orthogonal holomorphic bisectional curvature and
observed that the nonnegativity of the orthogonal holomorphic
bisectional curvature is preserved under the K$\ddot{a}$hler-Ricci
flow. (For the definition of the orthogonal holomorphic bisectional
curvature we will give in the following.) In 2006, X.X.Chen
\cite{ChenXX} generalized the Frankel conjecture in another aspect
with the orthogonal holomorphic bisectional curvature but under some
additional condition.  He proved that: \emph{any compact irreducible
K$\ddot{a}$hler manifold with positive orthogonal holomorphic
bisectional curvature and $c_1>0$ must be biholomorphic to the
complex projective space.}

\vskip 0.3cm \noindent {\bf Definition 1.1} \emph{ A complex
$n$-dimensional ($n\geq 2$) K$\ddot{a}$hler manifold $(M^n,h)$ is
said to have nonnegative orthogonal holomorphic bisectional
curvature if for any orthonormal basis $\{e_\alpha\}$, the following
holds:
$$R(e_\alpha,\overline{e_\alpha},e_\beta,\overline{e_\beta})
=R_{\alpha\bar{\alpha}\beta\bar{\beta}}\geq 0,$$ for any
$\alpha\neq\beta.$}

If we consider the K$\ddot{a}$hler manifold as a Riemannian
manifold, then we define the manifold has nonnegative orthogonal
holomorphic bisectional curvature by
$$R(u_i,Ju_i,Ju_j,u_j)\geq 0,$$ for any $<u_i,u_j>=<u_i,Ju_j>=0$,
where $J$ is the complex structure of $M$. The above Definition 1.1
is equivalent to that in the Riemannian case. Indeed, we can choose
an orthonormal basis $\{u_1,u_2,\cdots, u_{2n}\}$ such that
$Ju_i=u_{n+i}$ for $i=1,2,\cdots,n$. Set
$e_i=\frac{1}{\sqrt{2}}(u_i-\sqrt{-1}Ju_i)$, then $\{e_i\}$ is an
orthonormal basis. It follows that
$$R(e_i,\overline{e_i},e_j,\overline{e_j})=R_{i\bar{i}j \bar{j}}
=R(u_i,Ju_i,Ju_j,u_j),$$ for any $i\neq j$. This implies the two
definitions are equivalent.

Recently, Seshadri \cite{Se} gives the classification of manifolds
with nonnegative isotropic curvature. He proved that: \emph{any
compact irreducible K$\ddot{a}$hler manifold with nonnegative
isotropic curvature must be either a Hermitian symmetric manifold or
biholomorphic to the complex projective space.} From the computation
in Lemma 2.1 in \cite{Se}, we can see that nonnegative isotropic
curvature implies the nonnegative orthogonal holomorphic bisectional
curvature. However, the converse is not true. Following we give an
example and other examples can be given in a similar way:

\vskip 0.3cm \noindent {\bf Example 1.2}  Let
$$(M,h)=(\Sigma,g)\times (CP^n,g_0),$$ where $\Sigma$ is a Riemann
surface with Gauss curvature $\kappa(\Sigma)\geq-4$ and
$\min(\kappa(\Sigma))=-4$ and $g_0$ is the standard Fubini-Study
metric such that the sectional curvature of $CP^n$ satisfies $1\leq
K(p)\leq 4$. In the following, we want to show that $M$ has
nonnegative orthogonal holomorphic bisectional curvature but the
isotropic curvature is not nonnegative.

Indeed, suppose $\widehat{e_0}$ and $\{\widehat{e_i}\}, (1\leq i\leq
n),$ are the orthonormal basis of $T_p^{1,0}(\Sigma)$ and
$T_q^{1,0}(CP^n)$ respectively. Then we can naturally extend them to
be the orthonomal basis $\{e_i\}, (0\leq i\leq n),$ of
$T_x^{1,0}(M)$ at the point $x=(p,q)\in M$, such that $$pr_1
(e_0)=\widehat{e_0}, \mbox{ and }  pr_2(e_i)=\widehat{e_i},$$ where
$pr_1,pr_2$ denote the canonical projection onto $\Sigma$ and $CP^n$
respectively.

Now for any two orthogonal vectors $X,Y$ on $M$, we assume that:
$$X=\sum_{i=0}^na_ie_i, \mbox{ and
} Y=\sum_{i=0}^nb_ie_i,$$ where $a_i,b_i$ are complex numbers
satisfy $\sum\limits_{i=0}^na_i\overline{b_i}=0$.

Then by direct computation we can get that
$$\arraycolsep=1.5pt\begin{array}{rcl}
&&\hskip
0.5cmR(X,\overline{X},Y,\overline{Y})\\[4mm]
&&=|a_0|^2|b_0|^2R_{0\bar{0}0\bar{0}}
+4\sum\limits_{i=1}^n(|a_i|^2|b_i|^2)+2\sum\limits_{i=1}^n\sum\limits_{j\neq
i}
(|a_i|^2|b_j|^2+a_i\overline{a_j}b_j\overline{b_i})\\[4mm]
&&\geq
-4|\sum\limits_{i=1}^na_i\overline{b_i}|^2+4\sum\limits_{i=1}^n(|a_i|^2|b_i|^2)
+2\sum\limits_{i=1}^n\sum\limits_{j\neq i}
(|a_i|^2|b_j|^2+a_i\overline{a_j}b_j\overline{b_i})\\[4mm]
&&=2\sum\limits_{1\leq i<j\leq n}|a_ib_j-a_jb_i|^2\\[4mm]
&&\geq 0.
\end{array}$$
This implies that the orthogonal holomorphic bisectional curvature
of $M$ is nonnegative. On the other hand, by direct computation or
the result of \cite{MW}, it is easy to see that the isotropic
curvature is not nonnegative.

Clearly nonnegative holomorphic bisectional curvature also implies
the nonnegative orthogonal holomorphic bisectional curvature,
naturally we want to know the relations between holomorphic
bisectional curvature and isotropic curvature. By the work of Ivey
\cite{Ivey}, we know that in the complex 2-dimensional case,
nonnegative holomorphic bisectional curvature implies nonnegative
isotropic curvature. So the result of Seshadri \cite{Se} can be
viewed as a generalization of Mok's theorem on the generalized
Frankel conjecture in complex 2-dimension. But in higher dimensional
case, we do not know whether this is also true, since the
nonnegative holomorphic bisectional curvature means the bisectional
curvature is nonnegative on any holomorphic complex plane, while
nonnegative isotropic curvature requires on any 2-dimensional
isotropic plane. Even though, we know that both holomorphic
bisectional curvature and isotropic curvature imply orthogonal
holomorphic bisectional curvature. So orthogonal holomorphic
bisectional curvature is the weakest one among the three curvature
conditions. In \cite{ChenXX}, X.X.Chen asked a question: whether a
compact K$\ddot{a}$hler manifold with positive orthogonal
holomorphic bisectional curvature necessary has $c_1>0$. In this
paper, we give an affirmative answer to this question and hence
solve the Question/Conjecture 1.6 in \cite{ChenXX}. Moreover, we
will also give a complete classification of manifolds with
nonnegative orthogonal holomorphic bisectional curvature. This can
be considered as an extension of the generalized Frankel conjecture.
Our main result is the following:

\vskip 0.3cm \noindent {\bf Theorem 1.3} \emph{ Suppose $(M^n,h)$ is
an $n$-dimensional $(n\geq 2)$ compact K$\ddot{a}$hler manifold of
nonnegative orthogonal holomorphic bisectional curvature. Let
$(\tilde{M}^n,\tilde{h})$ be its universal covering space. Then
$(\tilde{M}^n,\tilde{h})$ is isometrically biholomorphic to one of
the following two cases:}

(1) \emph{ $(C^k,h_0)\times (M_1,h_1)\times \cdots \times
(M_l,h_l)\times (CP^{n_1},\theta_1)\times \cdots
\times(CP^{n_r},\theta_r)$, \vskip 0.1cm \noindent where $h_0$
denotes the Euclidean metric on $C^k, h_i (1\leq i\leq l)$ are
canonical metrics on the irreducible compact Hermitian symmetric
spaces $M_i$ of rank $\geq 2$, and $\theta_j (1\leq j\leq r)$ is a
K$\ddot{a}$hler metric on $CP^{n_j}$ carrying nonnegative orthogonal
holomorphic bisectional curvature;}

(2) \emph{ $(Y,g_0)\times (M_1,h_1)\times \cdots \times
(M_l,h_l)\times (CP^{n_1},\theta_1)\times \cdots
\times(CP^{n_r},\theta_r)$, \vskip 0.1cm \noindent where $Y$ is a
simply connected Riemann surface with Gauss curvature negative
somewhere or a simply connected noncompact K$\ddot{a}$hler manifold
with $\dim(Y)\geq 2$ and has nonnegative orthogonal holomorphic
bisectional curvature and the minimum of the holomorphic sectional
curvature $<0$ somewhere, $M_i, CP^{n_j} (1\leq i\leq l, 1\leq j\leq
r)$ are the same as in case (1). Moreover, we have the holomorphic
sectional curvatures of $M_i$ and $CP^{n_j}$ are $\geq-\min\{\mbox
{holomorphic sectional
 curvature of }$  $Y\}>0.$ }

This paper contains three sections and the organization is as
follows. In section 2, we will prove the positivity of the first
Chern class under the positive orthogonal holomorphic bisectional
curvature condition and give some results on the irreducible
manifolds which will be used in the proof of our main theorem. In
section 3, we will complete the proof of the Theorem 1.3.

\vskip 0.3cm \noindent {\bf Acknowledgement}  We would be indebted
to our advisor Professor X.P.Zhu for provoking our interest to
this problem. We are grateful to Professor B.L.Chen for many
suggestions and discussions. We would also like to thank Professor
S.H.Tang for discussions.

\begin{center}{ \bf \large \bf {2. Some Results on Irreducible Manifolds}}\end{center}

In the following we first give a similar result to \cite{MW} in
terms of the orthogonal holomorphic bisectional curvature in the
K$\ddot{a}$hler manifolds. We will show that the curvature term in
the Weitzenb$\ddot{o}$ck formula on $(1,1)$-forms involves only the
orthogonal holomorphic bisectional curvature. This also gives the
answer to the positivity of the first Chern class under the positive
orthogonal holomorphic bisectional curvature condition. In this
section we always assume that the complex dimension $n$ of the
K$\ddot{a}$hler manifold $M^n$ satisfies $n\geq 2$.

\vskip 0.1cm \noindent{\bf Theorem 2.1} \emph{ Let $(M^n,h)$ be a
compact K$\ddot{a}$hler manifold with nonnegative orthogonal
holomorphic bisectional curvature. Then all real harmonic
$(1,1)$-forms are parallel. Furthermore, we have}

(i) \emph{ If $b_{1,1}(M)=\dim H^{1,1}(M)=1$, then $c_1(M)>0$;}

(ii) \emph{ If in addition $M$ is locally irreducible, then we
have $ b_{1,1}(M)=\dim H^{1,1}(M)=1$ and hence by (i) we have
$c_1(M)>0$.}

\vskip 0.1cm \noindent{\bf Proof.} Suppose $(M^n,h)$ is a compact
K$\ddot{a}$hler manifold with nonnegative orthogonal holomorphic
bisectional curvature and $J$ is the complex structure. Let $\eta$
be a nontrivial harmonic $(1,1)$-form on $M$.

In the following, we want to show that $\eta$ is parallel. Indeed,
the parallelity of $\eta$ was already obtained by \cite{W} under the
condition of nonnegative holomorphic bisectional curvature. For the
completeness of our paper, we will adapt the argument in \cite{W}
and \cite{MW} to show that $\eta$ is parallel under the condition of
nonnegative orthogonal holomorphic bisectional curvature.

Now we can choose an orthonormal basis $\{e_\beta\}_{\beta=1}^n$
such that under this basis
$$\eta=\frac{\sqrt{-1}}{2}\sum\limits_{\beta=1}^n2a_\beta\cdot
e_\beta\wedge \overline{e_\beta}.$$ Set
$$e_\beta=\frac{1}{\sqrt{2}}(u_\beta-\sqrt{-1}Ju_\beta), \quad for \quad 1\leq\beta\leq n,$$
where $\{u_1,Ju_1,\cdots,u_n,Ju_n\}$  is an orthonormal basis of
$M$ in the sense of considering $M$ as a Riemannian manifold. So
in the basis $\{u_1,Ju_1,\cdots,u_n,Ju_n\}$, $\eta$ becomes
$$\eta=-\sum_{i=1}^na_i\cdot u_i\wedge Ju_i.$$
By the Bochner formula we have
$$\triangle \eta=\nabla^*\nabla\eta+\pounds(\eta),$$
where
$\pounds(\eta)=-\frac{1}{4}\sum\limits_{i,j}(R(\eta_i),\eta_j)
[\eta_i,[\eta_j,\eta]]$ and $(\cdot,\cdot)$ denotes the
corresponding Riemannian metric.  Then by the same argument as in
\cite{MW}, we know that
$$\arraycolsep=1.5pt\begin{array}{rcl}
&&(\pounds(\eta),\eta)=\frac{1}{2}\sum\limits_{\alpha>0}(R(X_\alpha),X_{-\alpha})
(-\alpha(\eta)X_{-\alpha},\alpha(\eta)X_\alpha)\\[4mm]
&&\hskip
2.1cm=\frac{1}{2}\sum\limits_{\alpha>0}-\alpha(\eta)^2(R(X_\alpha),X_{-\alpha}),
\end{array}$$
where $\alpha(\eta)$ satisfies
$[X_\alpha,\eta]=-[\eta,X_\alpha]=-\alpha(\eta)X_\alpha$ and
$-\alpha(\eta)^2$ is nonnegative and the symbols are the same as
in \cite{ChenXX}. Now for the positive roots $x_i+x_j, (1\leq
i<j\leq n),$ we have
$$ X_\alpha=\frac{1}{2}(u_i+\sqrt{-1}Ju_i)\wedge (u_j+\sqrt{-1}Ju_j)
=\overline{e_i}\wedge \overline{e_j},$$and
$$X_{-\alpha}=\frac{1}{2}(u_i-\sqrt{-1}Ju_i)\wedge (u_j-\sqrt{-1}Ju_j)
=e_i\wedge e_j.$$ For the positive roots $x_i-x_j, (1\leq i<j\leq
n),$ we have
 $$X_\alpha=\frac{1}{2}(u_i+\sqrt{-1}Ju_i)\wedge (u_j-\sqrt{-1}Ju_j)
=\overline{e_i}\wedge e_j,$$ and
$$X_{-\alpha}=\frac{1}{2}(u_i-\sqrt{-1}Ju_i)\wedge (u_j+\sqrt{-1}Ju_j)
=e_i\wedge \overline{e_j}.$$ So $(R(X_\alpha),X_{-\alpha})=0$ for
the previous case and for the other case, we have
$$(R(X_\alpha),X_{-\alpha})=R(\overline{e_i}\wedge
e_j,e_i\wedge\overline{e_j})
=R(\overline{e_i},e_j,\overline{e_j},e_i)=R(e_i,\overline{e_i},e_j,\overline{e_j})\geq
0,$$ since the orthogonal holomorphic bisectional curvature is
nonnegative. Then by the standard Bochner argument we can obtain
that all real harmonic $(1,1)$-forms are parallel.

In order to prove the left conclusions (i) and (ii), we evolve the
metric by the K$\ddot{a}$hler Ricci flow:
$$
      \left\{
       \begin{array}{lll}
\frac{\partial}{\partial t}g_{i\bar{j}}(x,t)=-R_{i\bar{j}}(x,t),
          \\[4mm]
  g_{i\bar{j}}(x,0)=h_{i\bar{j}}(x).
       \end{array}
    \right.
$$
Then by Shi's short-time existence theorem, we know that there is a
$T>0$ such that the Ricci flow has a smooth bounded curvature
solution $(M,g_{i\bar{j}}(t))$ for $t\in [0,T)$. It is due to
Cao-Hamilton \cite{CH} that the solution $g_{i\bar{j}}(t)$ still has
nonnegative orthogonal holomorphic bisectional curvature. Suppose
$\{e_\alpha\}$ is an orthonormal basis, then  for any
$\alpha\neq\beta$, we have:
$$\arraycolsep=1.5pt\begin{array}{rcl}
&&\hskip
0.5cmR(e_\alpha-e_\beta,\overline{e_\alpha}-\overline{e_\beta},
e_\alpha+e_\beta,\overline{e_\alpha}+\overline{e_\beta})\\[4mm]
&&=R_{\alpha\bar{\alpha}\alpha\bar{\alpha}}+R_{\beta\bar{\beta}\beta\bar{\beta}}
-R_{\alpha\bar{\beta}\alpha\bar{\beta}}-R_{\beta\bar{\alpha}\beta\bar{\alpha}}\\[4mm]
&&\geq 0,
\end{array}\eqno (2.1)$$ where we have used the assumption and the
result that the nonnegativity of the orthogonal holomorphic
bisectional curvature is preserved under the Ricci flow due to
Cao-Hamilton \cite{CH}. Similarly change $e_\beta$ by
$\sqrt{-1}e_\beta$, we have
$$R_{\alpha\bar{\alpha}\alpha\bar{\alpha}}+R_{\beta\bar{\beta}\beta\bar{\beta}}
+R_{\alpha\bar{\beta}\alpha\bar{\beta}}+R_{\beta\bar{\alpha}\beta\bar{\alpha}}\geq
0.\eqno (2.2)$$ By (2.1) and (2.2) we obtain that
$$R_{\alpha\bar{\alpha}\alpha\bar{\alpha}}+R_{\beta\bar{\beta}\beta\bar{\beta}}\geq 0\eqno
(2.3)$$ for any orthonormal 2-frames $\{e_\alpha,e_\beta\}$. So by
the assumption and (2.3), we have
$$R=\sum_{\alpha,\beta}R_{\alpha\bar{\alpha}\beta\bar{\beta}}
=\sum_\alpha\sum_{\beta\neq\alpha}R_{\alpha\bar{\alpha}\beta\bar{\beta}}
+\sum_\alpha R_{\alpha\bar{\alpha}\alpha\bar{\alpha}}\geq 0. \eqno
(2.4)$$

If $b_{1,1}(M)=\dim H^{1,1}(M)=1$, let $\rho$ and $\omega$ denote
the Ricci form and K$\ddot{a}$hler form respectively, then by the
Hodge theory, we have $\rho=\lambda\omega+\eta$, where $\lambda$
is a real number and $\int_M<\omega,\eta>=0$. On the other hand,
we have
$$\int_M<\rho,\omega>=\frac{1}{4}\int_M R=\lambda\parallel
\omega\parallel^2\geq 0,$$ since the scalar curvature $R\geq 0$ by
(2.4). Hence we have $c_1(M)\geq 0$. Moreover if the scalar
curvature $R$ at some point is positive, then $c_1(M)>0$. So now we
can assume that the scalar curvature $R(t)\equiv 0$ for all
sufficiently small $t$. Then by the evolution equation of the scalar
curvature
$$\frac{\partial R}{\partial t}=\triangle R+|Ric|^2$$
we know that for all sufficiently small $t$,
$$ Ric(t)\equiv 0. \eqno(2.5)$$

We claim that for all $\alpha$, the holomorphic sectional
curvature $R_{\alpha\bar{\alpha}\alpha\bar{\alpha}}=0.$

Indeed, by (2.3)-(2.5), we know that for any $\alpha\neq\beta$:
$$R_{\alpha\bar{\alpha}\beta\bar{\beta}}=0,\quad and \quad
R_{\alpha\bar{\alpha}\alpha\bar{\alpha}}+R_{\beta\bar{\beta}\beta\bar{\beta}}=0.$$
Suppose there exists $1\leq\alpha\leq n$ such that
$R_{\alpha\bar{\alpha}\alpha\bar{\alpha}}\neq 0$, then
$$R_{\alpha\bar{\alpha}}=\sum_{\beta\neq\alpha}R_{\alpha\bar{\alpha}\beta\bar{\beta}}
+R_{\alpha\bar{\alpha}\alpha\bar{\alpha}}\neq 0.$$ And this
contradicts with (2.5). So we have proved the claim and hence the
curvature operator is equal to zero. Therefore $(M^n,h)$ is flat.
However, note that $n\geq 2$, we know that there exists no compact
and flat K$\ddot{a}$hler manifold satisfying $b_{1,1}(M)=\dim
H^{1,1}(M)=1$. Thus the scalar curvature must be positive at some
point. Hence $c_1(M)>0$. This completes the proof of (i).

In the following we will give the proof of (ii). We argue by
contradiction. Suppose $b_{1,1}(M)=\dim H^{1,1}(M)>1$, then by the
same argument as in \cite{MW} in the proof of Theorem 2.1 (b) and
note that $M$ is locally irreducible, we know that $h$ is
hyper-K$\ddot{a}$hler and hence Ricci flat. So by the argument
above, we know that $M$ is flat. And this is a contradiction with
the local irreducibility of $M$. So $b_{1,1}(M)=\dim
H^{1,1}(M)=1$. Then by (i) we know that $c_1(M)>0$. This completes
the proof of (ii).

Therefore we complete the proof of Theorem 2.1.$$\eqno \#$$

From Theorem 2.1 and the result of \cite{ChenXX}, we immediately
obtain:

\vskip 0.1cm \noindent{\bf Corollary 2.2} \emph{ Let $(M^n,h)$ be a
compact K$\ddot{a}$hler manifold with positive orthogonal
holomorphic bisectional curvature. Then the first Chern class
$c_1(M)>0$. Moreover, the underlying manifold is  biholomorphic to
$CP^n$.} \vskip 0.1cm \noindent{\bf Proof.} Since $(M^n,h)$ has
positive orthogonal holomorphic bisectional curvature, we get that
$M$ is locally irreducible. Then by Theorem 2.1 (ii) we know that
$c_1(M)>0$. Combining the result of \cite{ChenXX}, we obtain that
$M$ is biholomorphic to the complex projective  space $CP^n$.
$$\eqno \#$$

Suppose $(M^n,h)$ is a compact K$\ddot{a}$hler manifold with
nonnegative orthogonal holomorphic bisectional curvature and
$g_{i\bar{j}}(t), 0\leq t\leq \delta,$ is the solution to the
K$\ddot{a}$hler Ricci flow with the initial data $h$. Let $P$ be the
bundle with the fixed metric $h$ and the fibre over $p\in M$
consists of all orthogonal 2-vectors $\{X,Y\}\subset T_p^{1,0}(M)$.
We define a function on $P\times (0,\delta)$ by
$$u(\{X,Y\},t)=R(X,\overline{X},Y,\overline{Y}),$$ where $R$
denotes the pull-back of the curvature tensor of
$g_{i\bar{j}}(t)$.
 \vskip 0.1cm \noindent{\bf Proposition 2.3} \emph{ There exists $c>0$ such that
 $$\frac{\partial u}{\partial t}\geq Lu+c\min\bigg\{0,\inf_{\xi\in V,|\xi|=1}D^2u(\xi,\xi)\bigg\}-
 c\sup_{\xi\in V,|\xi|=1}Du(\xi)-cu,$$ where $L$ is the horizontal Laplacian on $P$ and
$V$ denotes the vertical subspace of the bundle.} \vskip 0.1cm
\noindent{\bf Proof.} According to Hamilton \cite{Ha86}, under the
evolving orthonormal frame $\{e_\alpha\}$, we have
$$\arraycolsep=1.5pt\begin{array}{rcl}
&&\frac{\partial}{\partial t}R_{\alpha\bar{\alpha}\beta\bar{\beta}}=
\triangle R_{\alpha\bar{\alpha}\beta\bar{\beta}}+\sum\limits
_{\mu,\nu}\bigg(R_{\alpha\bar{\alpha}\mu\bar{\nu}}R_{\nu\bar{\mu}\beta\bar{\beta}}
-|R_{\alpha\bar{\mu}\beta\bar{\nu}}|^2+|R_{\alpha\bar{\beta}\mu\bar{\nu}}|^2\bigg)\\[4mm]
&&\hskip 1.6cm=\triangle
R_{\alpha\bar{\alpha}\beta\bar{\beta}}+\sum\limits
_{\mu,\nu=\alpha,\beta}\bigg(R_{\alpha\bar{\alpha}\mu\bar{\nu}}R_{\nu\bar{\mu}\beta\bar{\beta}}
-|R_{\alpha\bar{\mu}\beta\bar{\nu}}|^2+|R_{\alpha\bar{\beta}\mu\bar{\nu}}|^2\bigg)\\[4mm]
&&\hskip 1.6cm+\bigg(\sum\limits_{\mu=\alpha,\beta
\atop\nu\neq\alpha,\beta} +\sum\limits_{\nu=\alpha,\beta
\atop\mu\neq\alpha,\beta}\bigg)\bigg(R_{\alpha\bar{\alpha}\mu\bar{\nu}}R_{\nu\bar{\mu}\beta\bar{\beta}}
-|R_{\alpha\bar{\mu}\beta\bar{\nu}}|^2+|R_{\alpha\bar{\beta}\mu\bar{\nu}}|^2\bigg)\\[4mm]
&&\hskip 1.6cm+\sum\limits_{\mu,\nu\neq\alpha,\beta}
\bigg(R_{\alpha\bar{\alpha}\mu\bar{\nu}}R_{\nu\bar{\mu}\beta\bar{\beta}}
-|R_{\alpha\bar{\mu}\beta\bar{\nu}}|^2+|R_{\alpha\bar{\beta}\mu\bar{\nu}}|^2\bigg)\\[4mm]
&&\hskip 1.6cm \stackrel{\Delta}{=} \triangle
R_{\alpha\bar{\alpha}\beta\bar{\beta}}+(I)+(II)+(III).
\end{array}\eqno (2.6)$$

During the following proof, we assume that $c$ denotes the various
positive constants which depend on the bound of the curvature and
its derivatives.

\vskip 0.1cm \noindent{\bf Claim 1.} \emph{There exist constants
$c_1>0,c_2>0$ such that
$$I\geq -c_1\cdot u(\{e_\alpha,e_\beta\},t)
-c_2\sup_{\xi\in V,|\xi|=1}Du(\{e_\alpha,e_\beta\},t)(\xi).$$}
Indeed: by definition and direct computation, we have
$$\arraycolsep=1.5pt\begin{array}{rcl}
&&I=\sum\limits
_{\mu,\nu=\alpha,\beta}(R_{\alpha\bar{\alpha}\mu\bar{\nu}}R_{\nu\bar{\mu}\beta\bar{\beta}}
-|R_{\alpha\bar{\mu}\beta\bar{\nu}}|^2+|R_{\alpha\bar{\beta}\mu\bar{\nu}}|^2)\\[4mm]
&&\hskip
0.3cm=R_{\alpha\bar{\alpha}\beta\bar{\beta}}(R_{\alpha\bar{\alpha}\alpha\bar{\alpha}}+
R_{\beta\bar{\beta}\beta\bar{\beta}}-R_{\alpha\bar{\alpha}\beta\bar{\beta}})
+|R_{\alpha\bar{\beta}\alpha\bar{\beta}}|^2+2Re(R_{\alpha\bar{\alpha}\alpha\bar{\beta}}
\overline{R_{\alpha\bar{\beta}\beta\bar{\beta}}})
\end{array}\eqno (2.7)$$
Now we consider the orthogonal 2-frames $\{\cos\theta
e_\alpha+\sin\theta e_\beta,-\sin\theta e_\alpha+\cos\theta
e_\beta\}$, we have
$$\arraycolsep=1.5pt\begin{array}{rcl}
&&\hskip 0.5cm u(\{\cos\theta e_\alpha+\sin\theta
e_\beta,-\sin\theta e_\alpha+\cos\theta
e_\beta\},t)\\[4mm]
&&=R(\cos\theta e_\alpha+\sin\theta e_\beta,\cos\theta
\overline{e_\alpha}+\sin\theta \overline{e_\beta},-\sin\theta
e_\alpha+\cos\theta
e_\beta,-\sin\theta\overline{e_\alpha}+\cos\theta\overline{e_\beta}).
\end{array}$$
Then
$$\frac{du}{d\theta}\bigg|_{\theta=0}=2Re(R_{\alpha\bar{\beta}\beta\bar{\beta}}
-R_{\alpha\bar{\alpha}\alpha\bar{\beta}}).$$ So
$$|Re(R_{\alpha\bar{\beta}\beta\bar{\beta}}
-R_{\alpha\bar{\alpha}\alpha\bar{\beta}})|\leq c\cdot\sup_{\xi\in
V,|\xi|=1}Du(\{e_\alpha,e_\beta\},t)(\xi),\eqno (2.8)$$ for some
constant $c>0$.

Similarly, if we change $e_\beta$ by $\sqrt{-1}e_\beta$, and
consider the orthogonal 2-frames $\{\cos\theta e_\alpha+\sin\theta
\sqrt{-1}e_\beta,-\sin\theta e_\alpha+\cos\theta
\sqrt{-1}e_\beta\}$, note that
$$\arraycolsep=1.5pt\begin{array}{rcl}
&R(\cos\theta e_\alpha+\sin\theta \sqrt{-1}e_\beta,\cos\theta
\overline{e_\alpha}-\sin\theta
\sqrt{-1}\overline{e_\beta},\\[4mm]
&\hskip 1.5cm-\sin\theta e_\alpha+\cos\theta
\sqrt{-1}e_\beta,-\sin\theta\overline{e_\alpha}-\cos\theta\sqrt{-1}\overline{e_\beta})\\[4mm]
&\hskip 0.4cm=R(\cos\theta e_\alpha+\sin\theta
\sqrt{-1}e_\beta,\cos\theta \overline{e_\alpha}-\sin\theta
\sqrt{-1}\overline{e_\beta},\\[4mm]
&\hskip 1.5cm\sqrt{-1}\sin\theta e_\alpha+\cos\theta
e_\beta,-\sqrt{-1}\sin\theta\overline{e_\alpha}+\cos\theta\overline{e_\beta}),
\end{array}$$
we can obtain that
$$|Im(R_{\alpha\bar{\beta}\beta\bar{\beta}}
-R_{\alpha\bar{\alpha}\alpha\bar{\beta}})|\leq c\cdot\sup_{\xi\in
V,|\xi|=1}Du(\{e_\alpha,e_\beta\},t)(\xi).\eqno (2.9)$$ By (2.8)
and (2.9) we get that
$$|R_{\alpha\bar{\beta}\beta\bar{\beta}}
-R_{\alpha\bar{\alpha}\alpha\bar{\beta}}|^2\leq c\cdot\sup_{\xi\in
V,|\xi|=1}Du(\{e_\alpha,e_\beta\},t)(\xi)$$ i.e.,
$$|R_{\alpha\bar{\beta}\beta\bar{\beta}}|^2+|R_{\alpha\bar{\alpha}\alpha\bar{\beta}}|^2
-2Re(R_{\alpha\bar{\alpha}\alpha\bar{\beta}}\overline{R_{\alpha\bar{\beta}\beta\bar{\beta}}})
\leq c\cdot\sup_{\xi\in V,|\xi|=1}Du(\{e_\alpha,e_\beta\},t)(\xi).$$
So we have
$$2Re(R_{\alpha\bar{\alpha}\alpha\bar{\beta}}\overline{R_{\alpha\bar{\beta}\beta\bar{\beta}}})
\geq -c\cdot\sup_{\xi\in
V,|\xi|=1}Du(\{e_\alpha,e_\beta\},t)(\xi),\eqno(2.10)$$ for some
constant $c>0$.

By (2.7) and (2.10), we know that
$$I\geq -c_1\cdot u(\{e_\alpha,e_\beta\},t)
-c_2\sup_{\xi\in V,|\xi|=1}Du(\{e_\alpha,e_\beta\},t)(\xi),$$ for
some constants $c_1>0,c_2>0$. So we have proved Claim 1. \vskip
0.1cm \noindent{\bf Claim 2.} \emph{There exists constant $c_3>0,$
such that $$II\geq -c_3\sup_{\xi\in
V,|\xi|=1}Du(\{e_\alpha,e_\beta\},t)(\xi).$$} Indeed: by
definition and direct computation, we have
$$\arraycolsep=1.5pt\begin{array}{rcl}
&&II=\bigg(\sum\limits_{\mu\neq\alpha,\beta
\atop\nu=\alpha,\beta}+\sum\limits _{\nu\neq\alpha,\beta
\atop\mu=\alpha,\beta}\bigg)\bigg(R_{\alpha\bar{\alpha}\mu\bar{\nu}}R_{\nu\bar{\mu}\beta\bar{\beta}}
-|R_{\alpha\bar{\mu}\beta\bar{\nu}}|^2+|R_{\alpha\bar{\beta}\mu\bar{\nu}}|^2\bigg)\\[4mm]
&&\hskip
0.5cm=\sum\limits_{\mu\neq\alpha,\beta}2Re\bigg(R_{\alpha\bar{\alpha}\alpha\bar{\mu}}
R_{\mu\bar{\alpha}\beta\bar{\beta}}+R_{\alpha\bar{\alpha}\beta\bar{\mu}}
R_{\mu\bar{\beta}\beta\bar{\beta}}\bigg)\\[4mm]
&&\hskip
0.5cm+\sum\limits_{\mu\neq\alpha,\beta}\bigg(|R_{\alpha\bar{\beta}\mu\bar{\beta}}|^2
-|R_{\alpha\bar{\mu}\beta\bar{\beta}}|^2+|R_{\alpha\bar{\beta}\alpha\bar{\mu}}|^2
-|R_{\alpha\bar{\alpha}\beta\bar{\mu}}|^2\bigg).
\end{array}\eqno (2.11)$$
Now for $\mu\neq\alpha,\beta$, we consider the orthogonal
2-vectors $\{e_\alpha+se_\mu,e_\beta\}$, we have
$$u(\{e_\alpha+se_\mu,e_\beta\},t)=R(e_\alpha+se_\mu,
\overline{e_\alpha}+s\overline{e_\mu},e_\beta,\overline{e_\beta}).$$
Then
$$\frac{du}{ds}\bigg |_{s=0}=R_{\mu\bar{\alpha}\beta\bar{\beta}}
+R_{\alpha\bar{\mu}\beta\bar{\beta}}
=2Re(R_{\alpha\bar{\mu}\beta\bar{\beta}}).$$ So we have
$$|Re(R_{\alpha\bar{\mu}\beta\bar{\beta}})|\leq c\sup_{\xi\in
V,|\xi|=1}Du(\{e_\alpha,e_\beta\},t)(\xi),\eqno(2.12)$$ for some
constant $c>0$.

Change $e_\mu$ by $\sqrt{-1}e_\mu$, we  can obtain
$$|Im(R_{\alpha\bar{\mu}\beta\bar{\beta}})|\leq c\sup_{\xi\in
V,|\xi|=1}Du(\{e_\alpha,e_\beta\},t)(\xi).\eqno(2.13)$$ By (2.12)
and (2.13) we get $$|R_{\alpha\bar{\mu}\beta\bar{\beta}}|\leq
c\sup_{\xi\in
V,|\xi|=1}Du(\{e_\alpha,e_\beta\},t)(\xi).\eqno(2.14)$$ Similarly,
we can obtain that $$|R_{\alpha\bar{\alpha}\beta\bar{\mu}}|\leq
c\sup_{\xi\in
V,|\xi|=1}Du(\{e_\alpha,e_\beta\},t)(\xi).\eqno(2.15)$$ By (2.11),
(2.14) and (2.15) we know that $$II\geq -c\sup_{\xi\in
V,|\xi|=1}Du(\{e_\alpha,e_\beta\},t)(\xi),$$ for some constant
$c>0$. Hence we proved Claim 2. \vskip 0.1cm \noindent{\bf Claim
3.} \emph{There exists constant $c_4>0,$ such that $$III\geq
c_4\cdot\min\bigg\{0,\inf_{\xi\in
V,|\xi|=1}D^2u(\{e_\alpha,e_\beta\},t)(\xi,\xi)\bigg\}.$$} Indeed:
in the following we will prove that
$$\sum_{\mu,\nu\neq\alpha,\beta}\bigg(R_{\alpha\bar{\alpha}\mu\bar{\nu}}
R_{\nu\bar{\mu}\beta\bar{\beta}}-|R_{\alpha\bar{\mu}\beta\bar{\nu}}|^2\bigg)\geq
c_4\cdot\min\bigg\{0,\inf_{\xi\in
V,|\xi|=1}D^2u(\{e_\alpha,e_\beta\},t)(\xi,\xi)\bigg\},$$ for some
constant $c_4>0$.

For any vectors $\omega_\alpha,\omega_\beta$ orthogonal to
$e_\alpha,e_\beta$, we define an orthogonal 2-vectors
$\{v_\alpha(s),v_\beta(s)\}$ by:
$$v_\alpha(s)=e_\alpha+s\omega_\alpha-\frac{1}{2}s^2
\sum_{j=\alpha,\beta}<\omega_\alpha,\omega_j>e_j+O(s^3),$$
$$v_\beta(s)=e_\beta+s\omega_\beta-\frac{1}{2}s^2
\sum_{j=\alpha,\beta}<\omega_\beta,\omega_j>e_j+O(s^3).$$ Then
consider
$$u(\{v_\alpha(s),v_\beta(s)\},t)=R(v_\alpha,\overline{v_\alpha},v_\beta,\overline{v_\beta}).$$
By direct computation we have
$$\arraycolsep=1.5pt\begin{array}{rcl}
&&\frac{1}{2}\frac{d^2u(s)}{ds^2}|_{s=0}=
R(\omega_\alpha,\overline{\omega_\alpha},e_\beta,\overline{e_\beta})
+R(e_\alpha,\overline{e_\alpha},\omega_\beta,\overline{\omega_\beta})
+2Re(R(\omega_\alpha,\overline{e_\alpha},e_\beta,\overline{\omega_\beta}))\\[4mm]
&&\hskip 1.9cm
+2Re(R(e_\alpha,\overline{\omega_\alpha},e_\beta,\overline{\omega_\beta}))
-(<\omega_\alpha,\omega_\alpha>+<\omega_\beta,\omega_\beta>)
R_{\alpha\bar{\alpha}\beta\bar{\beta}}\\[4mm]
&&\hskip
1.9cm-Re(<\omega_\alpha,\omega_\beta>R_{\beta\bar{\alpha}\beta\bar{\beta}}
+<\omega_\beta,\omega_\alpha>R_{\alpha\bar{\alpha}\alpha\bar{\beta}}).
\end{array}$$
So we have
$$\arraycolsep=1.5pt\begin{array}{rcl}
&&\hskip
0.3cmR(\omega_\alpha,\overline{\omega_\alpha},e_\beta,\overline{e_\beta})
+R(e_\alpha,\overline{e_\alpha},\omega_\beta,\overline{\omega_\beta})
+2Re(R(\omega_\alpha,\overline{e_\alpha},e_\beta,\overline{\omega_\beta}))\\[4mm]
&&
+2Re(R(e_\alpha,\overline{\omega_\alpha},e_\beta,\overline{\omega_\beta}))
-(<\omega_\alpha,\omega_\alpha>+<\omega_\beta,\omega_\beta>)
R_{\alpha\bar{\alpha}\beta\bar{\beta}}\\[4mm]
&&-Re(<\omega_\alpha,\omega_\beta>R_{\beta\bar{\alpha}\beta\bar{\beta}}
+<\omega_\beta,\omega_\alpha>R_{\alpha\bar{\alpha}\alpha\bar{\beta}})\\[4mm]
&&\geq c\cdot\min\bigg\{0,\inf\limits_{\xi\in
V,|\xi|=1}D^2u(\{e_\alpha,e_\beta\},t)(\xi,\xi)\bigg\},
\end{array}\eqno(2.16)$$
for some constant $c>0$. If we change $e_\alpha$ by
$-\sqrt{-1}e_\alpha$ and $e_\beta$ by $\sqrt{-1}e_\beta$, we can
obtain that
$$\arraycolsep=1.5pt\begin{array}{rcl}
&&\hskip
0.3cmR(\omega_\alpha,\overline{\omega_\alpha},e_\beta,\overline{e_\beta})
+R(e_\alpha,\overline{e_\alpha},\omega_\beta,\overline{\omega_\beta})
-2Re(R(\omega_\alpha,\overline{e_\alpha},e_\beta,\overline{\omega_\beta}))\\[4mm]
&&
+2Re(R(e_\alpha,\overline{\omega_\alpha},e_\beta,\overline{\omega_\beta}))
-(<\omega_\alpha,\omega_\alpha>+<\omega_\beta,\omega_\beta>)
R_{\alpha\bar{\alpha}\beta\bar{\beta}}\\[4mm]
&&+Re(<\omega_\alpha,\omega_\beta>R_{\beta\bar{\alpha}\beta\bar{\beta}}
+<\omega_\beta,\omega_\alpha>R_{\alpha\bar{\alpha}\alpha\bar{\beta}})\\[4mm]
&&\geq c\cdot\min\bigg\{0,\inf\limits_{\xi\in
V,|\xi|=1}D^2u(\{e_\alpha,e_\beta\},t)(\xi,\xi)\bigg\},
\end{array}\eqno(2.17)$$
By (2.16) and (2.17) we have:
$$\arraycolsep=1.5pt\begin{array}{rcl}
&&\hskip
0.3cmR(\omega_\alpha,\overline{\omega_\alpha},e_\beta,\overline{e_\beta})
+R(e_\alpha,\overline{e_\alpha},\omega_\beta,\overline{\omega_\beta})
+2Re(R(e_\alpha,\overline{\omega_\alpha},e_\beta,\overline{\omega_\beta}))\\[4mm]
&&\geq c\cdot\min\bigg\{0,\inf\limits_{\xi\in
V,|\xi|=1}D^2u(\{e_\alpha,e_\beta\},t)(\xi,\xi)\bigg\}.
\end{array}\eqno(2.18)$$
If we set
$$\arraycolsep=1.5pt\begin{array}{rcl}
&&A(X,\overline{Y})=R(X,\overline{Y},e_\beta,\overline{e_\beta}),\\[4mm]
&&B(X,Y)=R(\overline{e_\alpha},X,\overline{e_\beta},Y),\\[4mm]
&&C(X,\overline{Y})=R(e_\alpha,\overline{e_\alpha},X,\overline{Y}).
\end{array}$$
Then by (2.18) we know that
$$
\left(
\begin{array}{cccc}A & \overline{B}\\
B^T & \overline{C}
\end{array}
\right)\geq c\cdot\min\bigg\{0,\inf_{\xi\in
V,|\xi|=1}D^2u(\{e_\alpha,e_\beta\},t)(\xi,\xi)\bigg\}.
$$
Hence we have
$$tr(AC)-tr(B\overline{B})\geq c\cdot\min\bigg\{0,\inf_{\xi\in
V,|\xi|=1}D^2u(\{e_\alpha,e_\beta\},t)(\xi,\xi)\bigg\},$$ where
$c>0$ is a constant depending on the bound of the curvature and
its derivatives. i.e.,
$$\sum_{\mu,\nu\neq\alpha,\beta}\bigg(R_{\alpha\bar{\alpha}\mu\bar{\nu}}
R_{\nu\bar{\mu}\beta\bar{\beta}}-|R_{\alpha\bar{\mu}\beta\bar{\nu}}|^2\bigg)\geq
c\cdot\min\bigg\{0,\inf_{\xi\in
V,|\xi|=1}D^2u(\{e_\alpha,e_\beta\},t)(\xi,\xi)\bigg\},\eqno(2.19)$$
for some constant $c>0$.

By the definition of $III$ and (2.19), we get
$$III\geq
c\cdot\min\bigg\{0,\inf_{\xi\in
V,|\xi|=1}D^2u(\{e_\alpha,e_\beta\},t)(\xi,\xi)\bigg\},$$ for some
constant $c>0$. Therefore we have proved Claim 3.

By (2.6), Claim 1, Claim 2 and Claim 3, we can get that
$$\frac{\partial u}{\partial t}\geq Lu+c\min\bigg\{0,\inf_{\xi\in V,|\xi|=1}D^2u(\xi,\xi)\bigg\}-
 c\sup_{\xi\in V,|\xi|=1}Du(\xi)-cu,$$ for some constant $c>0$, where $L$ is the horizontal Laplacian on $P$ and
$V$ denotes the vertical subspace of the bundle.

This completes the proof of Proposition 2.3.$$\eqno \#$$

\vskip 0.1cm \noindent{\bf Remark 2.4} In our proof, we have used
the result that the nonnegativity of the orthogonal holomorphic
bisectional curvature is preserved under the K$\ddot{a}$hler Ricci
flow, which is due to Cao-Hamilton \cite{CH} in an unpublished work.
However, we only used this result for the first term of (2.7) and
for obtaining of (2.18). So if we assume
$R_{\alpha\bar{\alpha}\beta\bar{\beta}}=0$, then the first term of
(2.7) is equal to zero and (2.18) is also true. Then combining
$Du(\{e_\alpha,e_\beta\},t)=0$ and $D^2u(\{e_\alpha,e_\beta\},t)\geq
0$, we can see that the argument of Proposition 2.3 has already
given a proof to this result. Also it is not hard to see that the
positivity of the orthogonal holomorphic bisectional curvature is
preserved under the K$\ddot{a}$hler Ricci flow.

In the following we will give a result on the irreducible compact
K$\ddot{a}$hler manifold with nonnegative orthogonal holomorphic
bisectional curvature.

\vskip 0.1cm \noindent{\bf Proposition 2.5} \emph{Let $(M^n,h)$ be a
compact irreducible K$\ddot{a}$hler manifold with nonnegative
orthogonal holomorphic bisectional curvature. Then either $M$ is
biholomorphic to the complex projective space or $(M,h)$ is
isometrically biholomorphic to an irreducible compact Hermitian
symmetric manifold of rank $\geq 2$.} \vskip 0.1cm \noindent{\bf
Proof.} Suppose $(M^n,h)$ is a compact irreducible K$\ddot{a}$hler
manifold with nonnegative orthogonal holomorphic bisectional
curvature, then by Theorem 2.1 (ii), we know that $c_1(M)>0$.

First we
 evolve the metric by the K$\ddot{a}$hler Ricci flow:
$$
      \left\{
       \begin{array}{lll}
\frac{\partial}{\partial t}g_{i\bar{j}}(x,t)=-R_{i\bar{j}}(x,t),
          \\[4mm]
  g_{i\bar{j}}(x,0)=h_{i\bar{j}}(x).
       \end{array}
    \right.
$$
According to Bando \cite{Ban}, we know that the evolved metric
$g_{i\bar{j}}(t), t\in (0,T)$, remains K$\ddot{a}$hler. Then by the
result due to Cao-Hamilton \cite{CH}, we know that for $t\in (0,T)$,
$g_{i\bar{j}}(t)$ has nonnegative orthogonal holomorphic bisectional
curvature. Moreover, according to Hamilton \cite{Ha86}, under the
evolving orthonormal frame $\{e_\alpha\}$, we have
$$\frac{\partial}{\partial t}R_{\alpha\bar{\alpha}\beta\bar{\beta}}=
\triangle R_{\alpha\bar{\alpha}\beta\bar{\beta}}+\sum\limits
_{\mu,\nu}\bigg(R_{\alpha\bar{\alpha}\mu\bar{\nu}}R_{\nu\bar{\mu}\beta\bar{\beta}}
-|R_{\alpha\bar{\mu}\beta\bar{\nu}}|^2+|R_{\alpha\bar{\beta}\mu\bar{\nu}}|^2\bigg).$$

Suppose $(M,h)$ is not locally symmetric. In the following, we
want to show that $M$ is biholomorphic to the complex projective
space $CP^n$.

Since the smooth limit of locally symmetric space is also locally
symmetric, we can obtain that there exists $\delta \in (0,T)$ such
that $(M,g_{i\bar{j}}(t))$ is not locally symmetric for $t \in
(0,\delta)$. Combining the K$\ddot{a}$hlerity of $g_{i\bar{j}}
(t)$ and Berger's holonomy theorem and note that $c_1(M)>0$, we
know that the holonomy group Hol$(g(t))=U(n)$.

As above, let $P$ be the fiber bundle with the fixed metric $h$
and the fiber over $p \in M$ consists of all orthogonal 2-vectors
$\{X,Y\}\subset T_p^{1,0}(M)$. We define a function $u$ on
$P\times(0,\delta)$ by
$$u(\{X,Y\},t)=R(X,\overline{X},Y,\overline{Y}),$$ where $R$
denotes the pull-back of the curvature tensor of $g_{i\bar{j}}(t)$.
Clearly we have $u\geq 0$, since $(M,g_{i\bar{j}}(t))$ has
nonnegative orthogonal holomorphic bisectional curvature due to
Cao-Hamilton \cite{CH}. Denote $F=\{(\{X,Y\},t)|u(\{X,Y\},t)=0,X\neq
0, Y\neq 0\}\subset P\times (0,\delta)$ consists of all pairs
$(\{X,Y\},t)$ such that $\{X,Y\}$ has zero orthogonal holomorphic
bisectional curvature with respect to $g_{i\bar{j}}(t)$. By
Proposition 2.3, we know that
$$\frac{\partial u}{\partial t}\geq Lu+c\min\bigg\{0,\inf_{\xi\in V,|\xi|=1}D^2u(\xi,\xi)\bigg\}-
 c\sup_{\xi\in V,|\xi|=1}Du(\xi)-cu,$$ for some constant $c>0$, where $L$ is the horizontal Laplacian on $P$ and
$V$ denotes the vertical subspace of the bundle. By Proposition 2
in \cite{BS2}, we know that the set
$$F=\bigg\{(\{X,Y\},t) | u(\{X,Y\},t)=0,X\neq 0, Y\neq 0\bigg\}\subset P\times
(0,\delta)$$ is invariant under parallel transport.

Next, by adapting the argument in \cite{Gu}, we claim that
$R_{\alpha\bar{\alpha}\beta\bar{\beta}}>0$ for all $t\in
(0,\delta)$ and all $\alpha\neq\beta$.

Indeed, suppose not. Then
$R_{\alpha\bar{\alpha}\beta\bar{\beta}}=0$ for some $t\in
(0,\delta)$ and some $\alpha\neq\beta$. Therefore
$$(\{e_\alpha,e_\beta\},t)\in F.$$
Combining $R_{\alpha\bar{\alpha}\beta\bar{\beta}}=0$ and the
computation for (2.7), (2.11) and (2.19) in Proposition 2.3, it is
not hard to obtain that:
$$
      \left\{
       \begin{array}{lll}
\sum\limits_{\mu,\nu}(R_{\alpha\bar{\alpha}\mu\bar{\nu}}R_{\nu\bar{\mu}\beta\bar{\beta}}
-|R_{\alpha\bar{\mu}\beta\bar{\nu}}|^2)=0,
          \\[4mm]
R_{\alpha\bar{\beta}\mu\bar{\nu}}=0, \quad \forall \mu, \nu,
          \\[4mm]
R_{\alpha\bar{\alpha}\mu\bar{\beta}}=
R_{\beta\bar{\beta}\mu\bar{\alpha}}=0,  \quad \forall \mu.
         \end{array}
    \right.
    \eqno(2.20)
$$
We define an orthonormal 2-frames
$\{\widetilde{e_\alpha},\widetilde{e_\beta}\}\subset T_p^{1,0}(M)$
by
$$\widetilde{e_\alpha}=\sin \theta\cdot e_\alpha-\cos \theta\cdot e_\beta,$$
$$\widetilde{e_\beta}=\cos \theta\cdot  e_\alpha+\sin \theta\cdot e_\beta.$$
Then
$$\overline{\widetilde{e_\alpha}}=
\sin \theta\cdot \overline{e_\alpha}-\cos \theta\cdot
\overline{e_\beta},$$
$$\overline{\widetilde{e_\beta}}=
\cos \theta\cdot \overline{e_\alpha}+\sin \theta\cdot
\overline{e_\beta}.$$ Since $F$ is invariant under parallel
transport and $(M,g_{i\bar{j}}(t))$ has holonomy group $U(n)$, we
obtain that
$$(\{\widetilde{e_\alpha},\widetilde{e_\beta}\},t)\in F,$$
that is,
$$R(\widetilde{e_\alpha},\overline{\widetilde{e_\alpha}},
\widetilde{e_\beta},\overline{\widetilde{e_\beta}})=0.$$ On the
other hand,
$$\arraycolsep=1.5pt\begin{array}{rcl}
&&R(\widetilde{e_\alpha},\overline{\widetilde{e_\alpha}},
\widetilde{e_\beta},\overline{\widetilde{e_\beta}})=\sin^2\theta\cos^2\theta
R_{\alpha\bar{\alpha}\alpha\bar{\alpha}}+\sin^3\theta\cos\theta
R_{\alpha\bar{\alpha}\alpha\bar{\beta}}+\sin^3\theta\cos\theta
R_{\alpha\bar{\alpha}\beta\bar{\alpha}}\\[4mm]
&&\hskip 3.3cm+\sin^4\theta
R_{\alpha\bar{\alpha}\beta\bar{\beta}}-\sin\theta\cos^3\theta
R_{\alpha\bar{\beta}\alpha\bar{\alpha}}-\sin^2\theta\cos^2\theta
R_{\alpha\bar{\beta}\alpha\bar{\beta}}\\[4mm]
&&\hskip 3.3cm-\sin^2\theta\cos^2\theta
R_{\alpha\bar{\beta}\beta\bar{\alpha}}-\sin^3\theta\cos\theta
R_{\alpha\bar{\beta}\beta\bar{\beta}}-\cos^3\theta\sin\theta
R_{\beta\bar{\alpha}\alpha\bar{\alpha}}\\[4mm]
&&\hskip 3.3cm-\sin^2\theta\cos^2\theta
R_{\beta\bar{\alpha}\alpha\bar{\beta}}-\sin^2\theta\cos^2\theta
R_{\beta\bar{\alpha}\beta\bar{\alpha}}-\cos\theta\sin^3\theta
R_{\beta\bar{\alpha}\beta\bar{\beta}}\\[4mm]
&&\hskip 3.3cm+\cos^4\theta
R_{\beta\bar{\beta}\alpha\bar{\alpha}}+\cos^3\theta\sin\theta
R_{\beta\bar{\beta}\alpha\bar{\beta}}+\cos^3\theta\sin\theta
R_{\beta\bar{\beta}\beta\bar{\alpha}}\\[4mm]
&&\hskip 3.3cm+\cos^2\theta\sin^2\theta
R_{\beta\bar{\beta}\beta\bar{\beta}}\\[4mm]
&&\hskip
2.9cm=\cos^2\theta\sin^2\theta(R_{\alpha\bar{\alpha}\alpha\bar{\alpha}}
+R_{\beta\bar{\beta}\beta\bar{\beta}}).
\end{array}$$
where in the last equality we have used (2.20). So we have
$$R_{\beta\bar{\beta}\beta\bar{\beta}}+R_{\alpha\bar{\alpha}\alpha\bar{\alpha}}=0,$$
if we choose $\theta$ such that $\cos^2\theta\sin^2\theta\neq 0$.

Clearly we can find an element of $U(n)$ such that it changes
$e_\alpha$ to $e_\mu$ and fixed $e_\beta$. Then we can see that
$$(\{e_\mu,e_\beta\},t)\in F.$$
By the same argument as above, we get
$$R_{\beta\bar{\beta}\mu\bar{\mu}}=R_{\beta\bar{\beta}\beta\bar{\beta}}+R_{\mu\bar{\mu}\mu\bar{\mu}}=0.$$
Similarly we can obtain that for any $e_\mu$ and $e_\nu$ with
$\mu\neq\nu$, the following holds:
$$R_{\mu\bar{\mu}\nu\bar{\nu}}=R_{\nu\bar{\nu}\nu\bar{\nu}}+R_{\mu\bar{\mu}\mu\bar{\mu}}=0.\eqno(2.21)$$
So we have the scalar curvature
$$R=\sum_{\alpha}R_{\alpha\bar{\alpha}}
=\sum_\alpha\sum_{\beta\neq\alpha}R_{\alpha\bar{\alpha}\beta\bar{\beta}}
+\sum_\alpha R_{\alpha\bar{\alpha}\alpha\bar{\alpha}}= 0. \eqno
(2.22)$$ Then by the same argument as in Theorem 2.1, we can obtain
that the manifold is flat, and this contradicts with the
irreducibility of $M$. Hence we prove that
$R_{\alpha\bar{\alpha}\beta\bar{\beta}}>0$, for all $t\in
(0,\delta)$ and all $\alpha\neq\beta$.

Therefore note that $c_1(M)>0$ and then using the result of
\cite{ChenXX}, we can get $M$ is biholomorphic to the complex
projective space $CP^n$.

This completes the proof of Proposition 2.5.$$\eqno \#$$

\begin{center}{ \bf \large \bf {3. The Proof of the Main Theorem}}\end{center}
\vskip 0.1cm \noindent{\bf Proof of Theorem 1.3.} Suppose $(M^n,h)$
is an $n$-dimensional $(n\geq 2)$ compact K$\ddot{a}$hler manifold
of nonnegative orthogonal holomorphic bisectional curvature. By
applying the standard de Rham decomposition theorem, we know that
the universal cover $(\tilde{M},\tilde{h})$ can be isometrically and
holomorphically splitted as $$(C^k,h_0)\times (M_1^{n_1},h_1)\times
\cdots \times (M_l^{n_l},h_l)$$ where each $(M_i^{n_i},h_i), 1 \leq
i \leq l ,$ is irreducible and non-flat, $h_0$ is the standard flat
metric on $C^k$ and $k, n_1,\cdots, n_l$ are nonnegative integers.

In the following we divide it into three cases:

\vskip 0.1cm \noindent{\bf Case 1.} \emph{ $k=0$ and in the de
Rham decomposition there exists a complex 1-dimensional
irreducible factor $\Sigma=M_1$ with Gauss curvature
$\kappa(\Sigma)$ negative somewhere}.

In this case, let $e_1$ be the unit basis of $T_p^{1,0}(\Sigma)$ and
$\{e_j^i\}, (1\leq j\leq n_i, 2\leq i\leq l)$, be the orthonormal
basis of $T_{q_i}^{1,0}(M_i)$ for arbitrary points $p\in \Sigma,
q_i\in M_i$. Naturally we can extend $e_1$ and $\{e_j^i\}$ to be an
orthonormal basis of $T^{1,0}_x(\tilde{M})$ for
$x=(p,q_2,\cdots,q_l)\in \tilde{M}$, still we denote by $e_1$ and
$e_j^i$,$(1\leq j\leq n_i, 2\leq i\leq l)$. Since $M$ has
nonnegative  orthogonal holomorphic bisectional curvature, we obtain
$$R(e_1-e_j^i,\overline{e_1}-\overline{e_j^i},e_1+e_j^i,\overline{e_1}+\overline{e_j^i})
=R^{(1)}_{1\bar{1}1\bar{1}}+R^{(i)}_{j\bar{j}j\bar{j}}=\kappa(p)+R^{(i)}_{j\bar{j}j\bar{j}}\geq
0,$$ where $R^{(i)}$ denotes the curvature on $M_i$. So for each
$i\neq 1$ we have
$$R^{(i)}_{j\bar{j}j\bar{j}}\geq -\kappa(p).$$ By the arbitrariness of
$p,q_i$, we know that
$$\min\{\mbox{holomorphic sectional  curvature of } M_i\} \geq
-\min\{\kappa(\Sigma)\}>0.$$ So we have proved that all $M_i, (i\neq
1),$ have nonnegative holomorphic bisectional curvature. If
$\dim(M_i)=n_i\geq 2$, then we know that it also has nonnegative
Ricci curvature. So $M_i$ is compact, otherwise, it will split off a
line and we can obtain a contradiction with the irreducibility of
$M_i$. Then by Proposition 2.5 we obtain that either $M_i$ is
biholomorphic to the complex projective space $CP^{n_i}$ or $M_i$ is
isometrically biholomorphic to an irreducible compact Hermitian
symmetric manifold of rank $\geq 2$. If $\dim(M_i)=n_i=1, (i\neq
1)$, then by the Gauss-Bonnet Theorem, we know that $M_i$ is
$S^2(=CP^1)$ with a nonnegatively curved metric. Hence this case is
contained in (2).

\vskip 0.1cm \noindent{\bf Case 2.} \emph{ $k=0$ and in the de
Rham decomposition there exists no complex 1-dimensional
irreducible factor or there exist complex 1-dimensional
irreducible factors and all these complex 1-dimensional
irreducible factors have nonnegatively curved metric}.

In this case, we know that all the complex 1-dimensional
irreducible factors, if exists, are compact by the Gauss-Bonnet
Theorem and are $S^2(=CP^1)$.

If all the irreducible factors $M_i$ with $\dim(M_i)\geq 2$ are
compact, then by Proposition 2.5 we obtain that either $M_i$ is
biholomorphic to the complex projective space $CP^{n_i}$ or $M_i$
is isometrically biholomorphic to an irreducible compact Hermitian
symmetric manifold of rank $\geq 2$. Hence this is contained in
(1).

If there exists an irreducible factor, without loss of generality,
denoted by $M_1$, is noncompact, then we claim that the minimal of
the holomorphic sectional curvature of $M_1$ $<0$ somewhere.
Otherwise, suppose the holomorphic bisectional curvature of $M_1$
$\geq 0$ and hence it has nonnegative Ricci curvature, so it is
compact which contradicts to the noncompactness of $M_1$. So we have
proved the claim. Then by the nonnegativity of the orthogonal
holomorphic bisectional curvature and the same argument as in Case
1, we know that all the other irreducible factors $M_i, (i\neq 1),$
have nonnegative holomorphic bisectional curvature and hence are
compact. Therefore as above, by Proposition 2.5 we obtain that
either $M_i, (i\neq 1),$ is biholomorphic to the complex projective
space $CP^{n_i}$ or $M_i, (i\neq 1),$ is isometrically biholomorphic
to an irreducible compact Hermitian symmetric manifold of rank $\geq
2$. This is contained in (2).

\vskip 0.1cm \noindent{\bf Case 3.} \emph{ $k\geq 1$}.

In this case, by the nonnegativity of the orthogonal holomorphic
bisectional curvature of $\tilde{M}$ and the same argument as in
Case 1, we know that all the other irreducible factors $M_i$ have
nonnegative holomorphic bisectional curvature. Again by the same
argument as in Case 1, we can obtain that if $\dim(M_i)=n_i\geq 2$,
then either $M_i$ is biholomorphic to the complex projective space
$CP^{n_i}$ or $M_i$ is isometrically biholomorphic to an irreducible
compact Hermitian symmetric manifold of rank $\geq 2$. If
$\dim(M_i)=n_i=1$, then by the Gauss-Bonnet Theorem, we know that
$M_i$ is $S^2(=CP^1)$ with a nonnegatively curved metric. This case
is contained in (1).

Hence from above argument, we have proved the Theorem 1.3. $$\eqno
\#$$

\end{document}